\documentclass[english,12pt]{article}
\usepackage{amsfonts}
\usepackage{amssymb}
\usepackage{amssymb, amsthm, amsmath}
\usepackage[latin1]{inputenc}
\usepackage{babel}
\usepackage{graphicx}
\usepackage{psfrag}
\usepackage{amstext}
\newtheorem{teorema}{\bf Theorem}

\newtheorem{corolario}{\bf Corollary}
\newtheorem{observacao}{\bf Remark}

\newtheorem{definicao}{ \bf Definition}

\newcommand{\abs}[1]{\lvert #1 \rvert}

\begin{document}
\title{A Congruence Theorem for Minimal Surfaces in $S^{5}$ with Constant Contact Angle}
\author{Rodrigo Ristow Montes and Jose A. Verderesi
\thanks{ristow@mat.ufpb.br and javerd@ime.usp.br}}
\date{Departamento de Matem\'atica , \\
      Universidade Federal da Para\'iba, \\[2mm]
      BR-- 58.051-900 ~ Jo\~ao Pessoa, P.B., Brazil \\
and \\
Departamento de Matem\'atica Pura, \\
Instituto de Matem\'atica e Estat\'{\i}stica, \\
Universidade de S\~ao Paulo, \\
Caixa Postal 66281, \\[2mm]
BR--05315-970~ S\~ao Paulo, S.P., Brazil}
\maketitle
\renewcommand{\thefootnote}{\fnsymbol{footnote}}
\renewcommand{\thefootnote}{\arabic{footnote}}
\setcounter{footnote}{0}
\thispagestyle{empty}
\begin{abstract}
\noindent
  We provide a congruence theorem for minimal surfaces in $S^5$ with constant contact angle using Gauss-Codazzi-Ricci equations. More precisely, we prove that Gauss-Codazzi-Ricci equations for minimal surfaces in $S^5$ with constant contact angle satisfy an equation for the Laplacian of the holomorphic angle. Also, we will give a characterization of flat minimal surfaces in $S^5$ with constant contact angle.
\end{abstract}

\smallskip
\noindent{\bf 2000 Math Subject Classification:} 53C42 - 53D10 - 53D35.
\section{Introduction}\mbox{}
We will establish a condition in order to investigate minimal surfaces in $S^5$ with constant contact angle using Gauss-Codazzi-Ricci equations. We define $\alpha$ to be the angle given by $ \cos \alpha = \langle ie_1 , v \rangle \nonumber$, where $e_1$ and $v$ are defined in section \ref{sec:section2}. The holomorphic angle $\alpha$ is the analogue of the K\"ahler angle introduced by Chern and Wolfson in \cite{CW} and \cite{W}.\\
In  \cite{MV} we introduced the notion of contact angle, that can be
considered as  a new geometric invariant useful to investigate
the geometry of immersed surfaces in $S^{3}$. Geometrically, the contact angle
$(\beta)$ is the complementary angle between the contact
distribution and the tangent space of the surface. Also in
\cite{MV}, we deduced formulas for the Gaussian curvature and the Laplacian
of an immersed minimal surface in $S^3$, and we gave a
characterization of the  Clifford Torus as the only minimal
surface in $S^3$
with constant contact angle.\\
Recently,  in \cite{RMV},  we construct a
family of minimal tori in $S^{5}$ with constant contact angle and
constant holomorphic angle. These tori are parametrized by the following circle
equation
\begin{equation}\label{eq:equacaoab}
a^2 + \left(b- \frac{\cos\beta}{1+\sin^2\beta}\right)^2 =  2\frac{\sin^4\beta}{(1+\sin^2\beta)^2},
\end{equation}
where  $a$ and $b$ are given in section \ref{parametro} (equation
(\ref{eq: segunda})).
In particular, when $a=0$ in \eqref{eq:equacaoab}, we recover the
examples found by Kenmotsu, in \cite{KK2}.  These examples are
defined for $0 < \beta < \frac{\pi}{2}$. Also, when $b=0$ in
\eqref{eq:equacaoab}, we find a new family of minimal tori in
$S^5$, and these tori are defined for $\frac{\pi}{4} < \beta <
\frac{\pi}{2}$. For $\beta = \frac{\pi}{2}$, we
give an alternative proof of this classification of a Theorem from 
Blair in \cite{Blair}, and Yamaguchi, Kon and Miyahara in
\cite{YKM} for Legendrian minimal surfaces in $S^5$ with
constant Gaussian curvature.\\
In this paper, we will establish a congruence theorem for  minimal surfaces in $S^5$ with constant contact angle ($\beta$). Using Gauss-Codazzi-Ricci equations, we prove the following theorem:
\begin{teorema}\label{teorema1}
Consider $S$ a riemannian surface, $\alpha : S \rightarrow ]0,\frac{\pi}{2}[$ a function over $S$ that verifies the following equation:
\begin{eqnarray}
\Delta(\alpha) &  =   & \cot\alpha\csc^3(\beta)\abs{\nabla\alpha}^2+a^2\cot\alpha\cot^4(\beta)-2a\cot\alpha\csc\beta\cot^2\beta\alpha_2\nonumber\\
&& -2\cos\alpha(\cot\beta-\tan\beta)\tan^2\beta\alpha_1+\sin\alpha\cos\alpha(5-\cot^4\beta-3\csc^2\beta)
\end{eqnarray}
then there exist one and  only one minimal immersion of $S$ into $S^5$ with constant contact angle such that $\alpha$ is the holomorphic angle of this immersion, where $a$ is given in section \ref{sec:section2}, and its determined as a function of $\alpha$ and $\beta$ in section \ref{sect:provas}.
\end{teorema}
As a immediatelly consequence of this method, we have a classification of flat minimal surfaces in $S^5$
\begin{corolario}
 Suppose that the contact angle ($\beta$) is constant and suppose that $S$ is a flat minimal surface in $S^5$ with constant principal curvature in the direction $e_3$,ie. $a$ is constant, then the holomorphic angle ($\alpha$) must be constant 
\end{corolario}
In general we have
\begin{observacao}
The Theorem \ref{teorema1} implies a more general classification in \cite{RMV}, in the sense that, in \cite{RMV}  we need the condition that both angles are constant.
\end{observacao}

\section{Contact Angle for Immersed Surfaces in  $S^{2n+1}$}\label{sec:section2}
Consider in $\mathbb{C}^{n+1}$ the following objects:
\begin{itemize}
\item the Hermitian product: $(z,w)=\sum_{j=0}^n z^j\bar{w}^j$;
\item the inner product: $\langle z,w \rangle = Re (z,w)$;
\item the unit sphere: $S^{2n+1}=\big\{z\in\mathbb{C}^{n+1} | (z,z)=1\big\}$;
\item the \emph{Reeb} vector field in $S^{2n+1}$, given by: $\xi(z)=iz$;
\item the contact distribution in $S^{2n+1}$, which is orthogonal to $\xi$:
\[\Delta_z=\big\{v\in T_zS^{2n+1} | \langle \xi , v \rangle = 0\big\}.\]
\end{itemize}
We observe that $\Delta$ is invariant by the complex structure of $\mathbb{C}^{n+1}$.

Let now $S$ be an immersed orientable surface in $S^{2n+1}$.
\begin{definicao}
The \emph{contact angle} $\beta$ is the complementary angle between the
contact distribution $\Delta$ and the tangent space $TS$ of the
surface.
\end{definicao}
Let $(e_1,e_2)$ be  a local frame of $TS$, where $e_1\in
TS\cap\Delta$, where $e_1$ is the characteristic field and  introduced by Bennequin in \cite{A}.\\
Then $\cos \beta = \langle \xi , e_2 \rangle
$. Finally, let $v$ be the unit vector in the direction of  the orthogonal projection  of $e_2$ on $\Delta$,
defined by the following relation
\begin{eqnarray}\label{eq:campoe2}
e_2 = \sin\beta v + \cos\beta \xi.
\end{eqnarray}
\section{Equations for Gaussian curvature and Laplacian of a minimal surface in
  $S^5$}\label{parametro}
In this section, we deduce the equations for the Gaussian curvature  and for the Laplacian
of a minimal surface in $S^5$ in terms of the
contact angle and the holomorphic angle.
Consider the normal vector fields
\begin{eqnarray}\label{eq:camposnormais}
e_3                   & = & i \csc \alpha e_1 - \cot \alpha v \nonumber\\
e_4                   & = & \cot \alpha e_1 + i \csc \alpha v \\
e_5                   & = & \csc \beta \xi - \cot \beta e_2
\nonumber
\end{eqnarray}
where $\beta \neq {{0,\pi}}$ and  $ \alpha \neq{{0,\pi}}$.
We will call $(e_j)_{1 \leq \ j \leq 5}$ an  \emph{adapted frame}.

Using (\ref{eq:campoe2}) and (\ref{eq:camposnormais}), we get
\begin{eqnarray}\label{eq:camposdistribution}\
v  =  \sin \beta e_2  - \cos \beta e_5 , \quad  iv   =  \sin \alpha e_4 - \cos \alpha e_1 \\
\xi  =  \cos \beta e_2 + \sin \beta e_5 \nonumber
\end{eqnarray}
It follows from (\ref{eq:camposnormais}) and
(\ref{eq:camposdistribution}) that
\begin{eqnarray}\label{eq:inormais}
ie_1 & = & \cos \alpha \sin \beta e_2 + \sin\alpha e_3 -\cos\alpha \cos \beta
e_5  \\
ie_2 & = & - \cos \beta z - \cos\alpha \sin\beta e_1  + \sin\alpha \sin
\beta e_4 \nonumber
\end{eqnarray}
Consider now the dual basis $(\theta^j)$ of $(e_j)$.
The connection forms $(\theta_k^j)$ are given by
\begin{eqnarray}
De_j    =   \theta_j^k e_k, \nonumber
\end{eqnarray}
and the second fundamental form  with respect to this frame are given by
\begin{equation}
\begin{array}{lclll}
II^j      & = & \theta_1^j \theta^1  +  \theta_2^j \theta^2; \quad
j=3, ..., 5 \nonumber.
\end{array}
\end{equation}
Using (\ref{eq:inormais}) and differentiating $v$ and $\xi$ on the
surface $S$, we get
\begin{eqnarray}\label{eq:dif}
D\xi & = & -\cos\alpha \sin\beta \theta^2 e_1 + \cos \alpha \sin \beta
 \theta^1 e_2 + \sin \alpha \theta^1 e_3 + \sin\alpha \sin \beta \theta^2 e_4\nonumber\\
          && - \cos \alpha \cos \beta \theta^1 e_5,  \\
Dv   & = & (\sin \beta \theta_2^1 - \cos \beta \theta_5^1)e_1 +
 \cos\beta(d\beta-\theta_5^2)e_2 + ( \sin\beta \theta_2^3 - \cos \beta
 \theta_5^3)e_3 \nonumber \\
          && +  ( \sin\beta \theta_4^2 - \cos \beta \theta_5^4)e_4 + \sin
 \beta(d\beta + \theta_2^5)e_5 \nonumber.
\end{eqnarray}
Differentiating $e_3$, $e_4$ and  $e_5$, we have
\begin{eqnarray}\label{eq:intrin}
\theta_3^1 & = &  -\theta_1^3 \nonumber \\
\theta_3^2 & = &  \phantom{-}\sin \beta ( d\alpha + \theta_4^1) - \cos \beta \sin \alpha
\theta^1 \nonumber \\
\theta_3^4 & = & \phantom{-} \csc \beta \theta_1^2 - \cot \alpha ( \theta_1^3 + \csc \beta
\theta_2^4) \nonumber \\
\theta_3^5 & = & \phantom{-} \cot \beta \theta_2^3 - \csc \beta \sin \alpha \theta^1 \nonumber \\
\theta_4^1 & = & -d\alpha - \csc \beta \theta_2^3 + \sin \alpha \cot \beta \theta^1 \nonumber \\
\theta_4^2 & = & - \theta_2^4 \nonumber \\
\theta_4^3 & = & \phantom{-} \csc \beta \theta_2^1 + \cot \alpha ( \theta_1^3 + \csc \beta
\theta_2^4) \nonumber \\
\theta_4^5 & = & \phantom{-} \cot \beta \theta_2^4 - \sin \alpha \theta^2 \nonumber \\
\theta_5^1 & = & -\cos \alpha \theta^2 -  \cot \beta \theta_2^1 \nonumber \\
\theta_5^2 & = & \phantom{-} d\beta + \cos \alpha \theta^1 \\  \label{eq: intrin}
\theta_5^3 & = & -\cot \beta \theta_2^3 + \csc \beta \sin \alpha \theta^1 \nonumber\\
\theta_5^4 & = & -\cot \beta \theta_2^4 + \sin \alpha \theta^2 \nonumber
\end{eqnarray}
The conditions of minimality and of symmetry  are equivalent to the following
equations:
\begin{eqnarray}\label{eq:minimal}
\theta_1^\lambda \wedge \theta^1 +  \theta_2^\lambda \wedge
\theta^2  = 0 = \theta_1^\lambda \wedge \theta^2 -
\theta_2^\lambda \wedge \theta^1.
\end{eqnarray}
On the surface $S$, we consider
\begin{eqnarray}
\theta_1^3 & = & a\theta^1 + b\theta^2 \nonumber
\end{eqnarray}
It follows from (\ref{eq:minimal}) that
\begin{eqnarray}\label{eq: segunda}
\theta_1^3 & = & \phantom{-} a\theta^1 + b\theta^2\nonumber\\
\theta_2^3 & = & \phantom{-} b\theta^1 - a\theta^2\nonumber\\
\theta_1^4 & = & \phantom{-} d\alpha + (b \csc\beta - \sin\alpha
\cot\beta)\theta^1 - a \csc \beta \theta^2\nonumber \\
\theta_2^4 & = & \phantom{-} d\alpha \circ J - a \csc\beta \theta^1 - (b
\csc\beta - \sin\alpha \cot\beta)\theta^2 \\
\theta_1^5   &  =  & \phantom{-} d\beta \circ J  - \cos \alpha \theta^2\nonumber\\
\theta_2^5   &  =  &  - d\beta - \cos \alpha \theta^1\nonumber
\end{eqnarray}
where $J$ is the complex structure  of $S$ is given by $Je_1=e_2$ and $Je_2=-e_1$.
Moreover, the normal connection forms are given by:
\begin{eqnarray}\label{eq:normalconexaobeta}
\theta_3^4 & = & - \sec\beta d\beta \circ J - \cot\alpha \csc\beta d\alpha
\circ J + a \cot\alpha \cot^2\beta \theta^1 \nonumber \\
&&+ ( b  \cot\alpha \cot^2\beta -
\cos\alpha \cot\beta \csc\beta + 2 \sec\beta \cos \alpha) \theta^2 \nonumber \\
\theta_3^5 & = & \phantom{-} (b \cot\beta - \csc\beta \sin\alpha) \theta^1 - a
\cot\beta \theta^2 \\
\theta_4^5 & = & \phantom{-} \cot\beta (d\alpha \circ J) - a \cot\beta
\csc\beta \theta^1 + ( -b\csc\beta \cot\beta + \sin\alpha (\cot^2\beta -1))\theta^2,\nonumber
\end{eqnarray}
while the Gauss equation is equivalent to the equation:
\begin{equation}\label{eq: Gauss}
\begin{array}{lcl}
d\theta_2^1 + \theta_k^1 \wedge \theta_2^k & = & \theta^1 \wedge \theta^2.
\end{array}
\end{equation}
Therefore, using  equations  (\ref{eq: segunda}) and (\ref{eq:
Gauss}), we have
\begin{eqnarray}\label{eq:curvatura1}
K  & = &  1 - |\nabla \beta|^2 - 2 \cos \alpha \beta_1 - \cos^2 \alpha
-(1+ \csc^2 \beta )(a^2 + b^2)
\nonumber\\
 &&+ 2b\sin\alpha\csc\beta\cot\beta  + 2 \sin \alpha \cot \beta \alpha_1 - |\nabla \alpha|^2 \nonumber\\
   && + 2 a \csc \beta
 \alpha_2 - 2b  \csc \beta \alpha_1 - \sin^2 \alpha \cot^2 \beta \nonumber\\
&&  =    1 - (1+ csc^2 \beta)(a^2+b^2)-2b\csc\beta
(\alpha_1-\sin\alpha \cot\beta)+2a \csc\beta \alpha_2
  \\
&& - |\nabla \beta + \cos \alpha e_1|^2 - |\nabla \alpha - \sin \alpha \cot
  \beta e_1|^2
        \nonumber
\end{eqnarray}
Using (\ref{eq:intrin}) and  the complex structure of $S$, we get
\begin{equation} \label{eq:conex}
\begin{array}{lcl}
\theta_2^1  &  =  &  \tan\beta(d\beta\circ J-2\cos\alpha\theta^2)
\end{array}
\end{equation}
Differentiating (\ref{eq:conex}), we conclude that
\begin{eqnarray}
d\theta_2^1 & = & ( -(1 + \tan^2
\beta)\abs{\nabla\beta}^2-\tan\beta\Delta\beta
-2\cos\alpha(1+2\tan^{2}\beta)\beta_1\nonumber\\
    &&  +2\tan\beta\sin\alpha\alpha_1-4\tan^{2}\beta\cos^{2}\alpha) \theta^1
    \wedge \theta^2 \nonumber
\end{eqnarray}
where $\Delta = tr \nabla^2 $ is the Laplacian of $S$. The Gaussian curvature is therefore
given by:
\begin{eqnarray}\label{eq:curvatura2}
K & = &  -(1 + \tan^2 \beta)\abs{\nabla\beta}^2-\tan\beta\Delta\beta
-2\cos\alpha(1+2\tan^{2}\beta)\beta_1 \nonumber\\
    &&
    +2\tan\beta\sin\alpha\alpha_1-4\tan^{2}\beta\cos^{2}\alpha.
\end{eqnarray}
From  (\ref{eq:curvatura1}) and (\ref{eq:curvatura2}), we obtain
the following formula for the Laplacian of $S$:
\begin{eqnarray}\label{eq:lapla}
\tan\beta\Delta\beta &  = & (1+ \csc^2 \beta)(a^2+b^2)+2b\csc\beta
(\alpha_1-\sin\alpha \cot\beta)-2a \csc\beta \alpha_2 \nonumber \\
&& - \tan^2 \beta( |\nabla \beta + 2 \cos \alpha e_1|^2
- |\cot \beta \nabla\alpha + \sin \alpha ( 1- \cot^2 \beta)e_1|^2 )\nonumber \\
 && + \sin^2 \alpha ( 1 -\tan^2 \beta)
\end{eqnarray}
\section{Gauss-Codazzi-Ricci equations for a minimal surface in $S^5$ with constant Contact angle $\beta$}\label{sect:equacoes}
In this section, we will compute Gauss-Codazzi-Ricci equations for a minimal surface in $S^{5}$ with constant Contact angle $\beta$.\\
Assume that the fundamental second form at the direction of the normal vector field $e_3$ diagonal. Then we have the following connection forms:
\begin{eqnarray}\label{eq: segunda1}
\theta_1^3 & = & \phantom{-} a\theta^1 \nonumber\\
\theta_2^3 & = & \phantom{-} - a\theta^2\nonumber\\
\theta_1^4 & = & \phantom{-} d\alpha + (- \sin\alpha
\cot\beta)\theta^1 - a \csc \beta \theta^2\nonumber \\
\theta_2^4 & = & \phantom{-} d\alpha \circ J - a \csc\beta \theta^1 - ( - \sin\alpha \cot\beta)\theta^2 \\
\theta_1^5   &  =  & \phantom{-} - \cos \alpha \theta^2\nonumber\\
\theta_2^5   &  =  & \phantom{-} - \cos \alpha \theta^1\nonumber
\end{eqnarray}
Normal connection forms are:
\begin{eqnarray}\label{eq:normalconexaobeta1}
\theta_3^4 & = & - \cot\alpha \csc\beta d\alpha
\circ J + a \cot\alpha \cot^2\beta \theta^1 \nonumber \\
&&+ ( - \cos\alpha \cot\beta \csc\beta + 2 \sec\beta \cos \alpha) \theta^2 \nonumber \\
\theta_3^5 & = & \phantom{-} (- \csc\beta \sin\alpha) \theta^1 - a
\cot\beta \theta^2 \\
\theta_4^5 & = & \phantom{-} \cot\beta (d\alpha \circ J) - a \cot\beta
\csc\beta \theta^1 + (  \sin\alpha (\cot^2\beta -1))\theta^2,\nonumber
\end{eqnarray}

Using the connection form (\ref{eq: segunda1}) and 
(\ref{eq:normalconexaobeta1}) in the Codazzi-Ricci equations, we have
\begin{eqnarray}
d\theta_1^3 + \theta_2^3 \wedge \theta_1^2 + \theta_4^3 \wedge \theta_1^4 + \theta_5^3 \wedge \theta_1^5  &  =  & 0 \nonumber
\end{eqnarray}
This implies that
\begin{eqnarray}\label{eq:Codazzi1}
&& -a_2 + a^2 (\cot\alpha \csc\beta \cot^2 \beta) - a \cot\alpha(\csc^2\beta + \cot^2\beta) \alpha_2 \\ 
&& - \cos\alpha \csc \beta (2(\cot\beta - \tan\beta)\alpha_1  - \sin\alpha (\cot^2 \beta -3)) + \cot\alpha \csc\beta |\nabla \alpha|^2 =  0\nonumber
\end{eqnarray}
Replacing the following (\ref{eq: segunda1}) and (\ref{eq:normalconexaobeta1}) in the Codazzi-Ricci equations
\begin{eqnarray}
d\theta_2^3 + \theta_1^3 \wedge \theta_2^1 + \theta_4^3 \wedge \theta_2^4 + \theta_5^3 \wedge \theta_2^5  &  =  & 0 \nonumber \\
d\theta_1^4 + \theta_2^4 \wedge \theta_1^2 + \theta_3^4 \wedge \theta_1^3 + \theta_5^4 \wedge \theta_1^5  &  =  & 0 \nonumber \\
d\theta_3^5 + \theta_1^5 \wedge \theta_3^1 + \theta_2^5 \wedge \theta_3^2 + \theta_4^5 \wedge \theta_3^4  &  =  & 0 \nonumber
\end{eqnarray}
We get
\begin{eqnarray}\label{eq:Codazzi2}
&& a_1 + a (\cot\alpha \alpha_1 + 6 \tan\beta
\cos\alpha)\nonumber \\
&& - 2 \sec\beta \cos\alpha \alpha_2=0
\end{eqnarray}
Using the connection form (\ref{eq: segunda1}) and (\ref{eq:normalconexaobeta1}) in the Codazzi-Ricci equations 
\begin{eqnarray}
d\theta_2^4 + \theta_1^4 \wedge \theta_2^1 + \theta_3^4 \wedge \theta_2^3 + \theta_5^4 \wedge \theta_2^5  &  =  & 0 \nonumber \\
d\theta_4^5 + \theta_1^5 \wedge \theta_4^1 + \theta_2^5 \wedge \theta_4^2 + \theta_3^5 \wedge \theta_4^3  &  =  & 0 \nonumber \\
d\theta_3^4 + \theta_1^4 \wedge \theta_3^1 + \theta_2^4 \wedge \theta_3^2 + \theta_5^4 \wedge \theta_3^5  &  =  & 0 \nonumber
\end{eqnarray}
We have
\begin{eqnarray}\label{eq:Codazzi3}
&& a_2 - a^2 (\cot\alpha \sin\beta \cot^2 \beta) + a \cot\alpha \alpha_2  \\
&& + 2 \cos\alpha (\cot\beta -3 \tan\beta)) + 2 \cos\alpha \sin\beta (\cot\beta -\tan\beta)\alpha_1 \nonumber\\
&& + \sin\alpha \cos\alpha \sin\beta (5 - \cot^2\beta) + \sin \beta \Delta \alpha  =  0\nonumber
\end{eqnarray}
Codazzi-Ricci equations
\begin{eqnarray}
d\theta_1^2 + \theta_3^2 \wedge \theta_1^3 + \theta_4^2 \wedge \theta_1^4 + \theta_5^2 \wedge \theta_1^5  &  =  & \theta^2 \wedge \theta^1 \nonumber \\
d\theta_1^5 + \theta_2^5 \wedge \theta_1^2 + \theta_3^5 \wedge \theta_1^3 + \theta_4^5 \wedge \theta_1^4  &  =  & 0 \nonumber
\end{eqnarray}
give the following equation
\begin{eqnarray}\label{eq:Codazzi4}
&& a^2 (1+\csc^2\beta) -2a\csc\beta\alpha_2  \nonumber\\
&&  + |\nabla \alpha|^2 + 2\sin\alpha (\tan \beta - \cot\beta) \alpha_1 - 4
\tan^2 \beta \cos^2 \alpha \nonumber \\
&& -\sin^2 \alpha (1-\cot^2\beta) = 0
\end{eqnarray}
The following Codazzi equation is automatically verified
\begin{eqnarray}
d\theta_2^5 + \theta_1^5 \wedge \theta_2^1 + \theta_3^5 \wedge \theta_2^3 + \theta_4^5 \wedge \theta_2^4  &  =  & 0 \nonumber
\end{eqnarray}
\section{Proof of the Results}\label{sect:provas}\mbox{}
In this section, we will prove the results using the equations of the sections 3 and 4
\subsection{Proof of Theorem 1}
Using equations (\ref{eq:Codazzi1}) and (\ref{eq:Codazzi3}) we have:
\begin{eqnarray}\label{eq:Codazzi6}
&& \sin^2\beta\Delta\alpha -\cot\alpha\csc\beta\abs{\nabla\alpha}^2-a^2\cot\alpha\cos^2\beta\cot^2\beta+2a\cot\alpha\sin\beta\cot^2\beta\alpha_2\nonumber\\ &&+2\cos\alpha(\cot\beta-\tan\beta)\cos^2\beta\alpha_1-\sin\alpha\cos\alpha(5\sin^2\beta-\cos^2\beta\cot^2\beta-3)=0 \nonumber
\end{eqnarray}
The equation above can be rewritten as the following equation
\begin{eqnarray}
\Delta(\alpha) &  =   & \cot\alpha\csc^3(\beta)\abs{\nabla\alpha}^2+a^2\cot\alpha\cot^4(\beta)-2a\cot\alpha\csc\beta\cot^2\beta\alpha_2\nonumber\\
&& -2\cos\alpha(\cot\beta-\tan\beta)\tan^2\beta\alpha_1+\sin\alpha\cos\alpha(5-\cot^4\beta-3\csc^2\beta)
\end{eqnarray}
Using the equation (\ref{eq:Codazzi4}) we get $a$ as a function of $\alpha$ and $\beta$
\begin{eqnarray}
&& a  = \frac{2\csc\beta\alpha_2+\sqrt f(\alpha,\beta)}{2(1+\csc^2\beta)}\nonumber
\end{eqnarray}
where $f(\alpha,\beta)$ is given by
\begin{eqnarray}
&&f(\alpha,\beta) = 4\alpha^2_2\cot^2\beta-4(1+\csc^2\beta)\alpha^2_1 \nonumber\\
&&-4(1+\csc^2\beta)(2\sin\alpha(\tan\beta-\cot\beta)\alpha_1-4\tan^2\beta\cos^2\alpha-\sin^2\alpha(1-\cot^2\beta))\nonumber
\end{eqnarray}
which prove Theorem 1.
\subsection{Proof of Corollary 1}
Suppose that $K=0$ at the equation ($\ref{eq:curvatura2}$), we determine $\alpha_1$ as a function of $\alpha$ 
\begin{eqnarray}\label{eq:alfa2}
 \alpha_1=2tg\beta\cot\alpha\cos\alpha
\end{eqnarray}
Now suppose that $a=k$, where $k$  is constant, at the equation ($\ref{eq:Codazzi2}$), we determine $\alpha_2$ as a function of $\alpha_1$
\begin{eqnarray}
 \alpha_2=k(\frac{\cos\beta\csc\alpha\alpha_1}{2}+3\sin\beta)
\end{eqnarray}
Using ($\ref{eq:alfa2}$) at the equation above, we have $\alpha_2$ as a function of $\alpha$
\begin{eqnarray}
\alpha_2=k\sin\beta(\cot^2\alpha+3)
\end{eqnarray}
Finally, using $\alpha_1$ and $\alpha_2$ at the equation ($\ref{eq:Codazzi1}$), we get $\alpha$ as a function of $\beta$, and therefore constant, which prove Corollary 1.

\end{document}